\documentclass[11pt, twoside, eqno]{article}
\usepackage{latexsym}
\usepackage{amssymb}
\usepackage{amsfonts}
\begin{document}
\begin{center}

\noindent {\bf \Large How we pass from fuzzy $po$-semigroups to fuzzy 
$po$-$\Gamma$-semigroups }\bigskip

\smallskip

{\bf Niovi Kehayopulu}\bigskip

{\small Department of Mathematics,
University of Athens \\
157 84 Panepistimiopolis, Athens, Greece \\
email: nkehayop@math.uoa.gr }
\end{center}
\date{ }

\bigskip

\noindent{\bf ABSTRACT.} The results on fuzzy ordered semigroups (or 
on fuzzy semigroups) can be transferred to fuzzy ordered gamma (or to 
fuzzy gamma) semigroups. We show the way we pass from fuzzy ordered 
semigroups to fuzzy ordered gamma semigroups.\bigskip

 \bigskip

\bigskip

The results of fuzzy ordered semigroups can be transferred into 
ordered gamma semigroups in the way indicated in the present paper. 
Let us give some results to justify what we say.

As we know, a fuzzy subset $f$ of an ordered groupoid $S$ is called a 
{\it fuzzy quasi-ideal} of $S$ if

(1) $(f\circ 1)\wedge (1\circ f)\preceq f$ and

(2) if $x\le y$, then $f(x)\ge f(y)$.\smallskip

While the quasi-ideals are defined via the multiplication and the 
order of fuzzy sets, the right, the left ideals and the bi-ideals are 
defined via the multiplication of $S$ as follows:

Let $(S,.,\le)$ be a $po$-groupoid. A fuzzy subset $f$ of $S$ is 
called a {\it fuzzy right ideal} of $S$ if

(1) $f(xy)\ge f(x)$ for every $x,y\in S$ and

(2) if $x\le y$, then $f(x)\ge f(y)$.

\noindent A fuzzy subset $f$ of $S$ is called a {\it fuzzy left 
ideal} of $S$ if

(1) $f(xy)\ge f(y)$ for every $x,y\in S$ and

(2) if $x\le y$, then $f(x)\ge f(y)$.

\noindent A fuzzy subset $f$ of an ordered semigroup $S$ is called a 
{\it bi-ideal} of $S$ if

(1) $f(xyz)\ge \min \{f(x),f(z)\}$ for every $x,y,z\in S$ and

(2) if $x\le y$, then $f(x)\ge f(y)$.

We have seen in [3] that a fuzzy subset of a groupoid $S$ is a fuzzy 
right (resp. left) ideal of $S$ if and only if the following 
assertions are satisfied:

(1) $f\circ 1\preceq f$ (resp. $1\circ f\preceq f$) and

(2) if $x\le y$, then $f(x)\ge f(y)$.

For a $po$-$\Gamma$-groupoid $(M,\Gamma,\le)$ we naturally have the 
following definitions: Any mapping $f$ from $M$ into the real closed 
interval [0,1] of real numbers is called a {\it fuzzy subset} of $M$. 
For a subset $A$ of $M$ the fuzzy subset $f_A$ is the characteristic 
function defined by:$$f_A : M\rightarrow [0,1] \mid x\rightarrow f_A 
(x):=\left\{ \begin{array}{l}
1\,\,\,\,\,\mbox { if }\,x\in A \\
0\,\,\,\,\,\,\mbox { if }\, x \not\in A.
\end{array} \right.$$For $A=\{a\}$ we write $f_a$ instead of 
$f_{\{a\}}$, and we have
$$f_a : M\rightarrow [0,1] \mid x\rightarrow f_a (x):=\left\{ 
\begin{array}{l}
1\,\,\,\,\,\mbox { if }\,x = a\\
0\,\,\,\,\,\,\mbox { if }\, x \ne a.
\end{array} \right.$$
For an element $a$ of $M$ we denote by $A_a$ the relation on $M$ 
defined by $$A_a:=\{(y,z) \mid a\le y\gamma z \mbox { for some } 
\gamma\in\Gamma\}.$$For two fuzzy subsets $f$ and $g$ of $M$ the 
multiplication $f\circ g$ is defined as follows: $$f\circ g : M 
\rightarrow [0,1] \mid a \rightarrow \left\{ \begin{array}{l}
\mathop  \bigvee \limits_{(y,z) \in {A_a}} \min\{f(y),g(z)\} 
\,\,\,\,\,\mbox { if }\,\,\,{A_a} \ne \emptyset \\
\;\;\;\;\;0\,\,\,\,\,\mbox { if }\,\,{A_a} = \emptyset \,
\end{array} \right.$$and in the set of all fuzzy subsets of $M$ we 
define the order relation as follows:$$f\preceq g \mbox { if and only 
if } f(a)\le g(a) \mbox { for all } a\in M.$$We also denote$$(f\wedge 
g)(a):=\min\{f(a),g(a)\} \mbox { for all } a\in M.$$
We denote by 1 the fuzzy subset of $M$ defined by: $1 : M \rightarrow 
[0,1] \mid x \rightarrow 1(x):=1$ and this is the greatest element of 
the set of fuzzy subsets of $M$. If $M$ is a $po$-$\Gamma$-groupoid 
(resp. $po$-$\Gamma$-semigroup), then the set of all fuzzy subsets of 
$M$ with the multiplication $``\circ"$ and the order $``\preceq"$ 
above is a $po$-$\Gamma$-groupoid (resp. $po$-$\Gamma$-semigroup).\\
For a subset $H$ of $M$ we denote by $(H]$ the subset of $M$ defined 
by: $$(H]:=\{t\in M \mid t\le a \mbox { for some } a\in H\}.$$
\noindent{\bf Definition 1.} Let $M$ be a $po$-$\Gamma$-groupoid.
 A fuzzy subset $f$ of $M$ is called a {\it fuzzy right ideal} of $M$ 
if

(1) $f(x\gamma y)\ge f(x)$ for every $x,y\in M$ and every 
$\gamma\in\Gamma$ and

(2) if $x\le y$, then $f(x)\ge f(y)$.\smallskip

\noindent{\bf Definition 2.} Let $M$ be a $po$-$\Gamma$-groupoid. A 
fuzzy subset $f$ of $M$ is called a {\it fuzzy left ideal} of $M$ if

(1) $f(x\gamma y)\ge f(y)$ for every $x,y\in M$ and every 
$\gamma\in\Gamma$ and

(2) if $x\le y$, then $f(x)\ge f(y)$. \smallskip

\noindent {\bf Definition 3.} A fuzzy subset $f$ of a 
$po$-$\Gamma$-semigroup $M$ is called a {\it fuzzy bi-ideal} of $S$ 
if

(1) $f(xyz)\ge \min \{f(x),f(z)\}$ for every $x,y,z\in M$ and

(2) if $x\le y$, then $f(x)\ge f(y)$.\smallskip

\noindent{\bf Proposition 4.} (cf. also [3]) {\it Let M be a 
$po$-$\Gamma$-groupoid. A fuzzy subset f of M is a fuzzy right ideal 
of M if and only if

$(1)$ $f\circ 1\preceq f$ and

$(2)$ if $x\le y$, then $f(x)\ge f(y)$.}\smallskip

\noindent{\bf Proof.} $\Longrightarrow$. Let $a\in S$. Then $(f\circ 
1)(a)\le f(a)$. In fact: If $A_a=\emptyset$, then $(f\circ 
1)(a):=0\le f(a)$. Let $A_a\not=\emptyset$. Then$$(f\circ 
1)(a):=\mathop  \bigvee \limits_{(y,z) \in {A_a}}\min\{f(y),1(z)\} 
=\mathop  \bigvee \limits_{(y,z) \in {A_a}}\min\{f(y)\}.$$On the 
other hand, $$f(y)\le f(a) \mbox { for every } (y,z)\in A_a.$$ 
Indeed: Let $(y,z)\in A_a$. Then $a\le y\gamma z$ for some 
$\gamma\in\Gamma$. Since $f$ is a fuzzy right ideal of $M$, we have 
$f(a)\ge f(y\gamma z)\ge f(y)$. Thus we have $$(f\circ 1)(a)=\mathop  
\bigvee \limits_{(y,z) \in {A_a}}\min\{f(y)\}\le 
f(a).$$$\Longleftarrow$. Let $x,y\in M$ and $\gamma\in\Gamma$. Then 
$f(x\gamma y)\ge f(x)$. Indeed: By hypothesis, we have $f(x\gamma 
y)\ge (f\circ 1)(x\gamma y)$. Since $(x,y)\in A_{x\gamma y}$, we have 
$$(f\circ 1)(x\gamma y):=\mathop  \bigvee \limits_{(u,v) \in 
{A_{x\gamma y}}}\min\{f(u),1(v)\}\ge \min\{f(x),1(y)\}=f(x).$$Thus 
$f(x\gamma y)\ge f(x)$.$\hfill\Box$\\In a similar way one can prove 
the following propositions:\smallskip

\noindent{\bf Proposition 5.} {\it Let M be a $po$-$\Gamma$-groupoid. 
A fuzzy subset f of M is a fuzzy left ideal of M if and only if

$(1)$ $1\circ f\preceq f$ and

$(2)$ if $x\le y$, then $f(x)\ge f(y)$.}\smallskip

\noindent{\bf Proposition 6.} {\it Let M be a $po$-$\Gamma$-groupoid. 
A fuzzy subset f of M is a fuzzy bi-ideal of M if and only if

$(1)$ $f\circ 1\circ f\preceq f$ and

$(2)$ if $x\le y$, then $f(x)\ge f(y)$.}\medskip

We characterize now the regular and the intra-regular 
$po$-$\Gamma$-semigroups in terms of fuzzy right and fuzzy left 
ideals.\smallskip

\noindent{\bf Definition 7.} A $po$-$\Gamma$-semigroup $M$ is called 
{\it regular} if $$a\in (a\Gamma M\Gamma a] \mbox { for every } a\in 
M.$$\noindent\\
\noindent{\bf Lemma 8.} {\it Let M be a regular 
$po$-$\Gamma$-semigroup. Then for every fuzzy right ideal f and every 
fuzzy subset g of M, we have $f\wedge g\preceq f\circ g$.}\smallskip

\noindent{\bf Proof.} Let $f$ be a fuzzy right ideal, $g$ a fuzzy 
subset of $M$ and $a\in M$. Since $M$ is regular, there exist $x\in 
M$ and $\gamma,\mu\in\Gamma$ such that $a\le (a\gamma x)\mu a$. Since 
$(a\gamma x,a)\in A_a$, we have $$(f\circ g)(a)=
\mathop  \bigvee \limits_{(y,z) \in {A_a}}\min\{f(y),g(z)\}\ge 
\min\{f(a\gamma x),g(a)\}.$$Since $f$ is a right ideal of $M$, we 
have $f(a\gamma x)\ge f(a)$. Then we have $$(f\circ g)(a)\ge 
\min\{f(a\gamma x),g(a)\}\ge \min\{f(a),g(a)\}=(f\wedge g)(a).$$
This holds for every $a\in M$, thus we have $f\wedge g\preceq f\circ 
g$.$\hfill\Box$\\
In a similar way we prove the following\smallskip

\noindent{\bf Lemma 9.} {\it Let M be a regular 
$po$-$\Gamma$-semigroup. Then for every fuzzy subset f and every 
fuzzy left ideal g of M, we have $f\wedge g\preceq f\circ 
g$.}\smallskip

\noindent{\bf Lemma 10.} {\it Let M be a $po$-$\Gamma$-groupoid, f a 
fuzzy right ideal and g a fuzzy left ideal of M. Then we have $f\circ 
g\preceq f\wedge g$.}\smallskip

\noindent{\bf Proof.} Let $a\in M$. Then $(f\circ g)(a)\le (f\wedge 
g)(a)$. In fact: If $A_a=\emptyset$, then $(f\circ g)(a):=0$. Since 
$a\in M$ and $f\wedge g$ is a fuzzy subset of $M$, we have $(f\wedge 
g)(a)\ge 0$, thus we have $(f\circ g)(a)\le (f\wedge g)(a)$. Let 
$A_a\not=\emptyset$. Then $$(f\circ g)(a):=\mathop  \bigvee 
\limits_{(y,z) \in {A_a}}\min\{f(y),g(z)\}.$$On the other hand,
$$\min\{f(y),g(z)\}\le (f\wedge g)(a) \;\;\forall\;\, (y,z)\in 
A_a.$$Indeed: Let $(y,z)\in A_a$. Then $y,z\in M$ and $a\le y\gamma 
z$ for some $\gamma\in\Gamma$. Since $f$ is a fuzzy right ideal of 
$M$, we have $f(a)\ge f(y\gamma z)\ge f(y)$. Since $g$ is a fuzzy 
left ideal of $M$, we have $g(a)\ge g(y\gamma z)\ge g(z)$. 
Thus$$\min\{f(y),g(z)\}\le \min\{f(a),g(a)\}=(f\wedge g)(a).$$Hence 
we obtain$$(f\circ g)(a)=\mathop  \bigvee \limits_{(y,z) \in 
{A_a}}\min\{f(y),g(z)\}\le (f\wedge g)(a),$$and the proof is 
complete.$\hfill\Box$\smallskip

\noindent{\bf Lemma 11.} (cf. also [1]) {\it Let M be a 
$po$-$\Gamma$-groupoid. Then A is a right (resp. left) ideal of M if 
and only if the characteristic function $f_A$ is a fuzzy right (resp. 
fuzzy left) ideal of M}.

For an element $a$ of $M$, we denote by $R(a)$, $L(a)$ the right and 
the left ideal of $M$ generated by $a$. We have $R(a)=(a\cup a\Gamma 
M]$ and $L(a)=(a\cup M\Gamma a]$.\smallskip

\noindent{\bf Lemma 12.} {\it A $po$-$\Gamma$-semigroup M is regular 
if and only if $$R(a)\cap L(a)\subseteq {\Big(}R(a)\Gamma L(a){\Big]} 
\mbox { for every } a\in M.$$ }
\noindent{\bf Theorem 13.} (see [2]) {\it A $po$-$\Gamma$-semigroup M 
is regular if and only if for every fuzzy right ideal f and every 
fuzzy left ideal g of M, we have$$f\wedge g\preceq f\circ g, \mbox { 
equivalently } f\wedge g=f\circ g.$$ }{\bf Proof.} Let $f$ be a fuzzy 
right ideal and $g$ a fuzzy left ideal of $M$. By Lemma 8, we have 
$f\wedge g\preceq f\circ g$. By Lemma 10, we have $f\circ g\preceq 
f\wedge g$. Then $f\wedge g=f\circ g$. Suppose $f\wedge g\preceq 
f\circ g$ for every right ideal $f$ and every left ideal $g$ of $M$. 
Then $M$ is regular. In fact: By Lemma 12, it is enough to prove that 
$R(a)\cap L(a)\subseteq (R(a)\Gamma L(a)]$ for every $a\in M$. Let 
$a\in M$ and $b\in R(a)\cap L(a)$. By Lemma 11, $f_{R(a)}$ is a fuzzy 
right ideal and $f_{L(a)}$ a fuzzy left ideal of $M$. By hypothesis, 
we have $$1=\min{\Big\{}f_{R(a)}(b),f_{L(a)}(b){\Big 
\}}={\Big(}f_{R(a)}\wedge f_{L(a)}{\Big)}(b)\le {\Big(}f_{R(a)}\circ 
f_{L(a)}{\Big)}(b),$$then $$1\le {\Big(}f_{R(a)}\circ 
f_{L(a)}{\Big)}(b).$$If $A_b =\emptyset$, then ${\Big(}f_{R(a)}\circ 
f_{L(a)}{\Big)}(b):=0$ which is impossible. If $A_b\not=\emptyset$,
then we have$${\Big (}f_{R(a)}\circ f_{L(a)}{\Big )}(b):=\mathop  
\bigvee \limits_{(y,z) \in {A_a}}\min\{f_{R(a)}(y),f_{L(a)}(z)\}.$$If 
$y\notin R(a)$ or $z\notin L(a)$ for every $(y,z)\in A_b$, 
then$$\min\{f_{R(a)}(y),f_{L(a)}(z)\}=0 \mbox { for every } (y,z)\in 
A_b,$$and then ${\Big (}f_{R(a)}\circ f_{L(a)}{\Big )}(b)=0$ which is 
impossible. Thus there exists $(y,z)\in A_b$ such that $y\in R(a)$ or 
$z\in L(a)$. And so, there exists $\gamma\in\Gamma$ such that $b\le 
y\gamma z\in R(a)\Gamma L(a)$, that is $b\in 
(R(a)\Gamma\L(a)]$.$\hfill\Box$
\smallskip

\noindent{\bf Definition 14.} A $po$-$\Gamma$-semigroup M is called 
{\it intra-regular} if $$a\in (M\Gamma a\Gamma a\Gamma M]\;\; 
\forall\, a\in M.$${\bf Lemma 15.} {\it A $po$-$\Gamma$-semigroup M 
is intra-regular if and only if$$R(a)\cap L(a)\subseteq 
{\Big(}L(a)\Gamma R(a){\Big]} \mbox { for every } a\in M.$$ }{\bf 
Theorem 16.} {\it A $po$-$\Gamma$-semigroup M is intra-regular if and 
only if for every fuzzy right ideal f and every fuzzy left ideal g of 
M, we have$$f\wedge g\preceq g\circ f.$$ }{\bf Proof.} 
$\Longrightarrow$. Let $f$ be a fuzzy right, $g$ a fuzzy left ideal 
of $M$ and $a\in M$. Since $M$ is intra-regular, there exist $x,y\in 
M$ and $\gamma,\mu,\rho\in\Gamma$ such that $a\le x\gamma a\mu a\rho 
y$. Since $(x\gamma a,a\rho y)\in A_a$, we have$$(g\circ 
f)(a):=\bigvee \limits_{(y,z) \in {A_a}}\min\{g(y),f(z)\}\ge \min 
\{g(x\gamma a),f(a\rho y)\}.$$Since $g$ is a fuzzy left ideal of $M$, 
$g(x\gamma a)\ge g(a)$. Since $f$ is a fuzzy right ideal of $M$, 
$f(a\rho y)\ge f(a)$. Thus we have$$(g\circ f)(a) \ge \min 
\{g(x\gamma a),f(a\rho y)\} \ge \min \{g(a),f(a)\}=(f\wedge 
g)(a),$$then $f\wedge g\preceq g\circ f$. \\$\Longleftarrow$. Let 
$a\in M$ and $b\in R(a)\cap L(a)$. By hypothesis, we 
have$$1=\min\{f_{R(a)}(b),f_{L(a)}(b)\}={\Big(}f_{R(a)}\wedge 
f_{L(b)}{\Big)}(b)\le {\Big(}f_{R(a)}\circ f_{R(b)}{\Big)}(b).$$Then 
$A_b\not=\emptyset$, and there exists $(y,z)\in A_b$ such that $y\in 
L(b)$ and $z\in R(b)$. Then there exists $\gamma\in\Gamma$ such that 
$b\le y\gamma z\in L(b)\Gamma R(b)$, so $b\in {\Big(}L(a)\Gamma 
R(a){\Big]}.$ By Lemma 15, $M$ is intra-regular.$\hfill\Box$

Let us characterize now the intra-regular, the regular and the left 
(right) regular $po$-$\Gamma$-semigroups in terms of fuzzy 
subsets.\smallskip

\noindent{\bf Lemma 17.} {\it Let M be a $po$-$\Gamma$-groupoid, f, g 
fuzzy subsets of M, and $a\in M$. The following are equivalent:

$(1)$ $(f\circ g)(a)\not=0$.

$(2)$ There exists $(x,y)\in A_a$ such that $f(x)\not=0$ and 
$g(y)\not=0$.}\smallskip

\noindent{\bf Proof.} $(1)\Longrightarrow (2)$. If $A_a=\emptyset$, 
then $(f\circ g)(a)=0$ which is impossible. So there exists $(x,y)\in 
A_a$, then$$(f\circ g)(a):=\mathop  \bigvee \limits_{(u,v) \in 
{A_{a}}} \min\{f(u),f(v)\}\ge \min\{f(x),g(y)\}.$$If $f(x)=0$ or 
$g(y)=0$, then $(f\circ g)(a)\ge 0$. Since $f\circ g$ is a fuzzy 
subset of $M$, we have $(f\circ g)(a)\le 0$. Then $(f\circ g)(a)=0$ 
which is impossible. Hence we have $f(x)\not=0$ and 
$g(y)\not=0$.$\hfill\Box$\smallskip

\noindent{\bf Corollary 18.} {\it Let M be a $po$-$\Gamma$-groupoid, 
f a fuzzy subset of M, and $a\in M$. The following are equivalent:

$(1)$ $(f\circ 1)(a)\not=0$.

$(2)$ There exists $(x,y)\in A_a$ such that $f(x)\not=0$.}\smallskip

\noindent{\bf Corollary 19.} {\it Let M be a $po$-$\Gamma$-groupoid, 
f a fuzzy subset of M, and $a\in M$. The following are equivalent:

$(1)$ $(1\circ g)(a)\not=0$.

$(2)$ There exists $(x,y)\in A_a$ such that $g(y)\not=0$.}\smallskip

\noindent{\bf Lemma 20.} {\it Let M be a $po$-$\Gamma$-groupoid and 
$a,b\in M$. Then we have$$b\le a\gamma a \mbox { for some } 
\gamma\in\Gamma\;\Longleftrightarrow\;(f_a\circ f_a)(b)\not=0.$$}{\bf 
Proof.} $\Longrightarrow$. Let $b\le a\gamma a$ for some 
$\gamma\in\Gamma$. Since $(a,a)\in A_b$, we 
have\begin{eqnarray*}(f_a\circ f_a)(b):=\mathop  \bigvee 
\limits_{(u,v) \in {A_{a}}} \min\{f_{a}(u),f_{a}(v)\}\ge 
\min\{f_{a}(a),f_{a}(a)\}=f_{a}(a)=1 \end{eqnarray*}
$\Longleftarrow$. Since $(f_a\circ f_a)(b)\not=0$, by Lemma 17, there 
exists $(x,y)\in A_b$ such that $f_a(x)\not=0$ and $f_a(y)\not=0$. 
Then $f_a(x)=f_a(y)=1$, and $x=y=a$. Since $b\le x\mu y$ for some 
$\mu\in\Gamma$, we have $b\le a\mu a$.$\hfill\Box$\smallskip

\noindent{\bf Lemma 21.} {\it Let M be a $po$-$\Gamma$-groupoid and f 
a fuzzy subset of M. Then$$f(a)\le (f\circ f)(a\gamma a) 
\;\;\forall\;a\in M\;\;\forall\;\gamma\in\Gamma.$$}
\noindent{\bf Proof.} Let $a\in M$ and $\gamma\in\Gamma$. Since 
$(a,a)\in A_a$, we have$$(f\circ f)(a\gamma a):=\mathop  \bigvee 
\limits_{(u,v) \in {A_{a\gamma a}}} \min\{f(u),f(v)\}\ge 
\min\{f(a),f(a)\}=f(a).$$Thus $f(a)\le (f\circ f)(a\gamma 
a)$.$\hfill\Box$\\We denote $f^2=f\circ f$.\smallskip

\noindent{\bf Lemma 22.} {\it Let M be a $po$-$\Gamma$-semigroup, f a 
fuzzy subset of M and $a\in M$. If $a\le x\mu a\gamma a\rho y$ for 
some $x,y\in M$, $\mu,\gamma,\rho\in\Gamma$, then $f(a)\le (1\circ 
f^2\circ 1)(a).$}\smallskip

\noindent{\bf Proof.} Since $(x\mu a\gamma a,y)\in A_a$, we have
\begin{eqnarray*}(1\circ f^2\circ 1)(a)&:=&\mathop  \bigvee 
\limits_{(u,v) \in {A_{a}}} \min\{(1\circ 
f^2)(u),1(v)\}\\&\ge&\min\{(1\circ f^2)(x\mu a\gamma 
a),1(y)\}\\&=&(1\circ f^2)(x\mu a\gamma a).\end{eqnarray*}Since 
$(x,a\gamma a)\in A_{x\mu a\gamma a}$, we 
have\begin{eqnarray*}(1\circ f^2)(x\mu a\gamma a)&:=&\mathop  \bigvee 
\limits_{(w,t) \in {A_{x\mu a\gamma a}}} 
\min\{1(w),f^2(t)\}\\&\ge&\min\{1(x),f^2(a\gamma a)\}\\&=&f^2(a\gamma 
a).\end{eqnarray*}Then, by Lemma 21, we have $(1\circ f^2\circ 
1)(a)\ge f^2(a\gamma a)\ge f(a)$.$\hfill\Box$
\smallskip

\noindent{\bf Theorem 23.} {\it A $po$-$\Gamma$-semigroup M is 
regular if and only if for every fuzzy subset $f$ of $M$, we have  
$f\preceq f\circ 1\circ f$.}\smallskip

\noindent{\bf Proof.} $\Longrightarrow$. Let $f$ be a fuzzy subset of 
$M$ and $a\in M$. Since $M$ is regular, there exist $x\in M$ and 
$\gamma,\mu\in\Gamma$ such that $a\le a\gamma x\mu a$. Since 
$(a\gamma x,a)\in A_a$, we have\begin{eqnarray*}(f\circ 1\circ 
f)(a)&:=&\mathop  \bigvee \limits_{(y,z) \in {A_a}}\min\{(f\circ 
1)(y),f(z)\}\\&\ge&\min\{(f\circ 1)(a\gamma 
x),f(a)\}.\end{eqnarray*}Since $(a,x)\in A_{a\gamma x}$, we 
have$$(f\circ 1)(a\gamma x):=\mathop  \bigvee \limits_{(u,v) \in 
{A_{a\gamma x}}}\min\{f(u),1(v)\}\ge 
\min\{f(a),1(x)\}=f(a).$$Then$$(f\circ 1\circ f)(a)\ge 
\{\min\{f(a),f(a)\}=f(a),$$so  $f\preceq f\circ 1\circ f$.\\
$\Longleftarrow$. let $a\in M$. Since $f_a$ is a fuzzy subset of $M$, 
by hypothesis, we have $a=f_a(a)\le (f_a\circ 1\circ f_a)(a)\le 1$, 
so ${\Big(}(f_a\circ 1)\circ f_a{\Big)}(a)\not=0$. By Lemma 17, there 
exists $(x,y)\in A_a$ such that $(f_a\circ 1)(x)\not=0$ and 
$f_a(y)\not=0$. By Corollary 18, there exists $(u,v)\in A_x$ such 
that $f_a(x)\not=0$. Since $a\le x\gamma y$ for some 
$\gamma\in\Gamma$ and $x=y=a$, we have $a\le a\gamma a\le (a\gamma 
a)\gamma a\in a\Gamma M\Gamma a$. Then $a\in (a\Gamma M\Gamma a]$, 
and $M$ is regular. $\hfill\Box$\smallskip

\noindent{\bf Theorem 24.} (see [4]) {\it A $po$-$\Gamma$-semigroup M 
is intra-regular if and only if, for every fuzzy subset f of M, we 
have$$f\preceq 1\circ f^2\circ 1.$$}{\bf Proof.} $\Longrightarrow$. 
Let $f$ be a fuzzy subset of $M$ and $a\in M$. Since $M$ is 
intra-regular, there exist $x,y\in M$ and $\mu,\gamma,\rho\in\Gamma$ 
such that $a\le x\mu a\gamma a\rho y$. Then, by Lemma 22, we have
$f(a)\le (1\circ f^2\circ 1)(a)$, so $f\preceq 1\circ f^2\circ 1$.\\
$\Longleftarrow$. Let $a\in M$. Since $f_a$ is a fuzzy subset of $M$, 
by hypothesis, we have $$1=f_a(a)\le (1\circ f^2_a\circ 1)(a)\le 1,$$ 
hence ${\Big(}(1\circ f^2_a)\circ 1{\Big)}(a)\not=0$. By Corollary 
18, there exists $(x,y)\in A_a$ such that $(1\circ f^2)(x)\not=0$.
By Corollary 19, there exists $(u,v)\in A_x$ such that 
$f_a^2(v)\not=0$. By Lemma 20, there exists $\gamma\in\Gamma$ such 
that $v\le a\gamma a$. In addition, $a\le x\mu y$ and $x\le u\rho v$ 
for some $\mu,\rho\in\Gamma$. Then we have $$a\le x\mu y\le (u\rho 
v)\mu y\le u\rho (a\gamma a)\mu y\in M\Gamma a\Gamma a\Gamma M,$$ and 
$a\in (M\Gamma a\Gamma a\Gamma M]$, thus $M$ is 
intra-regular.$\hfill\Box$\smallskip

\noindent{\bf Definition 25.} A $po$-$\Gamma$-semigroup $M$ is called 
{\it right regular} if $$a\in (a\Gamma a\Gamma M] \;\;\,\forall\;a\in 
M.$${\bf Theorem 26.} {\it A $po$-$\Gamma$-semigroup M is right 
regular if and only if, for every fuzzy subset f of M, we 
have$$f\preceq f^2\circ 1.$$}{\bf Proof.} $\Longrightarrow$.
Let $a\in M$. Since $M$ is right regular, there exist $x\in M$ and 
$\gamma,\mu\in\Gamma$ such that $a\le a\gamma a\mu x$. Since 
$(a\gamma a,x)\in A_a$, we have\begin{eqnarray*}(f^2\circ 
1)(a)&:=&\mathop  \bigvee \limits_{(u,v) \in {A_{a}}} \min\{ 
f^2(u),1(v)\}\\&\ge&\min\{f^2(a\gamma a),1(x)\}\\&=&f^2(a\gamma 
a).\end{eqnarray*}Since $(a,a)\in A_{a\gamma a}$, we have
$$f^2(a\gamma a)=\mathop  \bigvee \limits_{(w,t) \in {A_{a\gamma a}}} 
\min\{ f(w),f(t)\}\ge \min\{f(a),f(a)\}\\=f(a).$$So $(f^2\circ 
1)(a)\ge f(a)$, then $f\preceq f^2\circ 1$.\\
$\Longleftarrow$. Let $a\in M$. By hypothesis, we have $1=f_a(a)\le 
(f_a^2\circ 1)(a)$, so $(f_a^2\circ 1)(a)=1$. By Corollary 18, there 
exists $(x,y)\in A_a$ such that $f_a^2(x)\not=0$. By Lemma 20, there 
exists $\lambda\in\Gamma$ such that $x\le a\lambda a\le (a\lambda 
a)\lambda a\in a\Gamma a\Gamma M$, so $a\in (a\Gamma a\Gamma M]$.$
\hfill\Box$\smallskip

\noindent{\bf Definition 27.} A $po$-$\Gamma$-semigroup $M$ is called 
{\it left regular} if $$a\in (M\Gamma a\Gamma a] \;\;\,\forall\;a\in 
M.$${\bf Theorem 28.} {\it A $po$-$\Gamma$-semigroup M is left 
regular if and only if, for every fuzzy subset f of M, we 
have$$f\preceq 1\circ f^2.$$}

\noindent{\bf Proposition 29.} {\it If M is a $po$-$\Gamma$-groupoid 
and f a fuzzy right (resp. fuzzy left) ideal of M, then $f\circ 
1\preceq f$ (resp. $1\circ f\preceq f )$.}\smallskip

\noindent{\bf Proof.} Let $f$ be a fuzzy right ideal of $M$ and $a\in 
M$. Then $(f\circ 1)(a)\le f(a)$. Indeed: If $A_a=\emptyset$, then 
$(f\circ 1)(a):=0\le f(a)$. Let $A_a\not=\emptyset$. Then$$(f\circ 
1)(a):=\mathop  \bigvee \limits_{(x,y) \in {A_a}}\min\{f(x),1(y)\}=
\mathop  \bigvee \limits_{(x,y) \in {A_a}}\min\{f(x)\}.$$On the other 
hand, $f(x)\le f(a)$ for every $(x,y)\in A_a$. Indeed: Let $(x,y)\in 
A_a$. Then $a\le x\gamma y$ for some $\gamma\in\Gamma$. Since $f$ be 
a fuzzy right ideal of $M$, we have $$f(a)\ge f(x\gamma y)\ge f(x). 
$$Then $f\circ 1\preceq f$.$\hfill\Box$\smallskip

\noindent{\bf Corollary 30.} {\it If M is a $po$-$\Gamma$-groupoid, 
then the fuzzy right (and the fuzzy left) ideals of $M$ are 
subidempotent.} \smallskip

\noindent{\bf Proof.} Let $f$ be a fuzzy right ideal of $M$. Then 
$f^2\preceq f\circ 1\preceq f$.$\hfill\Box$.\smallskip

\noindent{\bf Corollary 31.} {\it If M is a regular 
$po$-$\Gamma$-semigroup, then the fuzzy right (and the fuzzy left) 
ideals of $M$ are idempotent.}\smallskip

\noindent{\bf Proof.}$\Longrightarrow$. Let $f$ be a fuzzy right 
ideal of $M$. Since $M$ is regular, by Theorem 23, we have $f\preceq 
(f\circ 1)\circ f\preceq f^2$. By Corollary 30, $f^2\preceq f$, thus 
we have $f=f^2$.$\hfill\Box$

 \bigskip

\noindent This paper under the title 
``On fuzzy $po$-$\Gamma$-semigroups" has been submitted in Armenian 
Journal of Mathematics on June 23, 2014 at 18:34 (the date in 
Greece).

\end{document}